\documentclass[11pt] {amsart}
\usepackage{amsmath,amsfonts,amssymb,amsxtra,latexsym,amscd,enumerate,amsthm}

\usepackage{latexsym}
\usepackage{epsfig}
\usepackage{verbatim}

\def\r{\mathcal{R}}

\def\a{\alpha}
\def\v{\varphi}

\def\d{\delta}
\def\b{\beta}

\def\p{\Phi}
\def\f{\varphi}
\def\ep{\epsilon}
\def\es{\emptyset}
\def\l{\lambda}
\def\o{\omega}

\def\R{\mathbb{R}}
\def\r{\mathcal{R}}

\def\I{\mathcal{I}}
\def\Q{\mathcal{Q}}
\def\Z{\mathbb{Z}}
\def\K{\mathcal{K}}

\def\P{\mathbb{P}}
\def\p{\mathcal{P}}

\def\cN{\mathcal{N}}

\def\S{\mathcal{S}}

\def\A{\mathcal{A}}
\def\B{\mathcal{B}}
\def\TT{\mathbb{T}}
\def\N{\mathbb{N}}
\def\Z{\mathbb{Z}}

\def\beq{\begin{equation}}
\def\eeq{\end{equation}}

\def\beq{\begin{equation}}
\def\eeq{\end{equation}}

\newtheorem{thm}{Theorem}

\newtheorem{obs}[thm]{Observation}

\newtheorem{d0}[thm]{Definition}

\newtheorem{l1}[thm]{Lemma}

\newtheorem{r1}[thm]{Remark}

\begin{document}
\title[Polynomial Carleson Operator in higher dimensions]{A Note on the Polynomial Carleson Operator in higher dimensions}

\author{Victor Lie}

\date{\today}

\address{Victor Lie, Department of Mathematics, Purdue, IN 46907 USA}

\thanks{The author was supported by the National Science Foundation under Grant No. DMS-1500958.}

\email{vlie@purdue.edu}

\address{Institute of Mathematics of the Romanian Academy, Bucharest, RO
70700 \newline \indent  P.O. Box 1-764}



\maketitle

%\hfill\parbox[l]{1.6 in}{\scriptsize \it{Dedicated to Elias Stein on the occasion of his $80^{th}$ birthday celebration.}}
%$\newline$
\begin{abstract}
We prove the $L^p$-boundedness, $1<p<\infty$, of the Polynomial Carleson operator in general dimension. This follows the author's resolution of the one dimensional case (\cite{q}, \cite{lvPC}) as well as the work of Zorin-Kranich (\cite{zk}) on the higher dimensional case in the setting $2\leq p<\infty$.

The techniques used in this paper are direct adaptations and natural extensions to the higher dimensional case of the one-dimensional methods developed in \cite{lvPC}.
\end{abstract}

\section{\bf Introduction}

In this paper we present the required adaptations of the main techniques developed in \cite{lvPC} in order to fully answer a conjecture of E. Stein regarding the boundedness properties of the so-called Polynomial Carleson operator:

$\newline$
{\bf Conjecture (\cite{s2},\cite{sw}).}
{\it Let $G$ denote either $\TT$ or $\R$ with $G^m:=\prod_{j=1}^m G$, $m\in\N$. Further, let $\Q_{d,m}$ be the class of all real-coefficient polynomials in $m$ variables with no constant term and of degree less than or equal to $d$, $d\in\N$, and let $K$ be a suitable Calder\'on--Zygmund operator on $G^m$. Then the Polynomial Carleson operator defined as
\beq\label{polcarl} C_{d,m}f(x):=\sup_{Q\in\Q_{d,m}}\left|
\,\int_{G^m}e^{i\,Q(y)}\,K(y)\,f(x-y)\,dy\,\right|\:
\eeq
obeys the bound
\beq\label{polcarlb}
\|C_{d,m} f\|_{L^p(G^m)} \lesssim \|f\|_{L^p(G^m)}
\eeq
for any $1<p<\infty$.}
$\newline$

For the motivation and history of the problem the interested reader is invited to consult the Introduction in \cite{lvPC}
as well as the references \cite{c1}, \cite{f}, \cite{hu}, \cite{lt3}, \cite{s2} and \cite{sw}.

For concreteness and simplicity assume\footnote{For more on this see Observation \ref{genkz} below. Also notice that if instead of $G=\R$ one chooses $G=\TT$ then $K$ on $\TT^m$ may be essentially regarded as the restriction to the unit cube centered at the origin of a Calderon-Zygmund kernel over $\R^m$.} in what follows that $K:\:\R^m\,\rightarrow\,\R$ is a Calderon-Zygmund kernel that obeys the following properties:
\beq\label{CZker}
\eeq
\begin{itemize}
\item $K$ is a tempered distribution which agrees with a $C^1$ function $K(x)$ for $x\not=0$;

\item the Fourier transform $\hat{K}$ is an $L^{\infty}$ function;

\item $|\partial_x ^{\a} K(x)|\leq C\,|x|^{-n+|\a|}$ for $0\leq |\a|\leq 1$ where here $C>0$ is an absolute constant.
\end{itemize}

With these notations and specifications the main result of our paper is:

$\newline$
{\bf Main Theorem.} {\it The above conjecture holds.}
$\newline$

\begin{obs}\label{genkz} The assumptions in \eqref{CZker} on the Calderon-Zygmund kernel could be significantly relaxed. Also, it is quite likely that one can extend our result to cover anisotropic Calderon-Zygmund kernels - for more on these see the Introduction in \cite{zk}. However, all these extensions seem to stress rather technical aspects than genuinely new ideas.
\end{obs}

The above theorem extends the previous results of the author in the one-dimensional case (see \cite{q} and \cite{lvPC}) and completes - with one key correction\footnote{For more on this see 3) in the Remarks section.} - the work of Zorin-Kranich, \cite{zk}, on the higher dimensional case in the $2\leq p<\infty$ regime.

As already mentioned earlier, the proof of our theorem follows very closely the corresponding approach in \cite{lvPC}. This is why we will insist here only on the more relevant adaptations needed and then simply recall/outline our proof in \cite{lvPC}.

The higher dimensional upgrade requires at least two such adaptations:

\begin{enumerate}
\item a proper generalization (to higher dimensions) of the time-frequency tile discretization of the Carleson operator;

\item a non-stationary phase principle for polynomial type phases - here used in the form of Van der Corput estimates as developed by E. Stein and S. Wainger in \cite{sw}.
\end{enumerate}

It is worth mentioning in this context that both (1) and (2) were exploited in \cite{zk}. As we will soon see, there are multiple - but essentially equivalent - ways to realize (1). \footnote{For a more elaborate discussion on this please see 2.1) in the Remarks section.}

For transparency and better understanding of the parallelism between the one and higher dimensional cases we choose to preserve the structure from our paper \cite{lvPC}.

\section{\bf Notations and construction of the tiles}\label{Not}

In this section we would like to make transparent the map/dictionary between the one-dimensional case and the corresponding higher dimensional case passing partly through the adaptations offered by Zorin-Kranich in \cite{zk}.

Consider the standard dyadic grids relative to the origin in $\TT^m$ and $\R^m$ respectively. A spacial dyadic cube $I$ in $\TT^m$ at scale $j\in\N$ has the form $I=\prod_{l=1}^{m}[2^{-j}\,k_l, 2^{-j}\,(k_l+1))$ with $\{k_l\}_{l=1}^{m}\subseteq\{0,\ldots\,2^j-1\}$ while a frequency dyadic cube at scale $-j$ has the form $\o=\prod_{l=1}^{m}[2^{j}\,r_l, 2^{j}\,(r_l+1))$ with $\{r_l\}_{l=1}^{m}\subseteq\Z$.

While in space we will decompose our universe in regular geometric shapes given by dyadic cubes in frequency the discretization will be realized in shapes that arise as intersections of tube-neighborhoods around suitable graph of polynomials. Regarded as part of the dual time-frequency representation once we fix a scale $j\in\N$ this double spacial/frequency discretization gives rise to a partition in \textit{tiles}.

To make the discussion more transparent, we structure our discussion as follows:

\subsection{The one dimensional case}

Throughout this section we consider $m=1$.  In this situation - maintaining the notations and definitions from \cite{lvPC} and abusing the language - the tiles $P\equiv\hat{P}$ are represented as neighborhoods\footnote{See in this direction the explicit interpretation provided by Observation 2 in \cite{lvPC}.} of volume one around the graphs of polynomials $q=\frac{dQ}{dx}\in \Q_{d-1,1}$ with $Q\in \Q_{d,1}$. Here, we performed an identification between the above geometric shape interpretation and the $(d+1)$-tuple of dyadic
intervals $P=[\a^1,\a^2,\ldots , \a^d,I]$ with $|\a^j|=|I|^{-1},\:\:j\in\{1,\ldots , d\}$ via the following convention:
setting
\beq\label{convbel00}
\Q_{d-1,1}(P):=\{q\in \Q_{d-1,1}\:|\:q(x^j_{I})\in \a^j\:\:\forall\:j\in\{1,\ldots,d\} \}\:,
\eeq
as (roughly) the $|I|^{-1}-$neighborhood of the \emph{central polynomial}
\beq\label{convbel0}
q_P(y):=\sum_{j=1}^{d}\frac{\prod_{k=1\atop{k\not=j}}^{d}(y-x_I^k)}{\prod_{k=1\atop{k\not=j}}^{d}(x_I^j-x_I^k)}\:c(\a^j)\:,
\eeq
and defining
\beq\label{convbel}
 q\in P\:\:\:\textrm{iff}\:\:\:q\in \Q_{d-1,1}(P)\,.
\eeq
one identifies $P$ with
\beq\label{gtile} \hat{P}:=\{(x,q(x))\:|\:x\in I\:\:
\&\:\: q\in P \}\:. \eeq

Equivalently, recalling from \cite{lvPC} the definition of the \emph{geometric factor} of the pair $(q,I)$ with $q=\frac{dQ}{dx}\in \Q_{d-1,1}$ and $I\subseteq\TT$ interval, that is
\beq\label{dqj}
\Delta_{Q'}(I)=\Delta_q(I):=|I|\,\|q\|_{L^{\infty}(I)}\:,
\eeq
we notice that the previous tile discretization is essentially equivalent with the following:

\textsf{Algorithm, $m=1$.}

\begin{itemize}
\item fix a scale $k\in\N$ and consider the set $\I(k)$ of all dyadic intervals $I\subseteq \TT$ at scale $k$, that is $|I|=2^{-k}$.

\item fix $I\in\I(k)$ and consider a maximal separated set of polynomials $\Q_{d-1,1}(I)\subset \Q_{d-1,1}$ such that
for any $q_{1},q_{2}\in\Q_{d-1,1}(I)$ one has
\beq\label{dqja}
\Delta_{q_1-q_2}(I)\geq 1\:.
\eeq

\item construct the uncertainty regions $\{P(q,I)\}_{q\in \Q_{d-1,1}(I)}$ partitioning the set of polynomials $\Q_{d-1,1}$ such that
  \beq\label{dqjaa0}
    P(q,I)\subset\{\tilde{q}\in \Q_{d-1,1}\,|\,\Delta_{\tilde{q}-q}(I)\leq \frac{1}{2}\}\:.
 \eeq
Notice that $P(q,I)$ corresponds to the set $\Q_{d-1,1}(P)$ introduced in \eqref{convbel00} with the central polynomial $q_P$ in \eqref{convbel0} corresponding to $q$ in  \eqref{dqjaa0}.

\item now we set a tile $P$ be a tuple of the form $(I, P(q,I))$ where $I$ ranges through the set of all dyadic intervals inside the torus while $q\in \Q_{d-1,1}(I)$. When more clarity is desired in order to associate $q,\,I$ with $P$ we will write $P=(I_P, P(q_P,I_P))$. The set of all tiles is denoted with $\P$.
\end{itemize}

\begin{r1}\label{remk} In the algorithm above we borrowed part of the language introduced in \cite{zk} by Zorin-Kranich in order to make the correspondence between the two approaches in \cite{lvPC} and \cite{zk} more transparent.
\end{r1}

\subsection{The general dimensional case}

For general $m$, the tile discretization could follow either of the above (morally equivalent) strategies with the obvious adaptations. However, it seems more advantageous and cleaner to follow the latter variant due to the less appealing form of Lagrange interpolation formulas for higher degree polynomials. This is the route that we embrace below:

We recall that in this case we deal with dyadic cubes $I\subseteq\TT^m$. From here on we follow the most natural adaptation of what we've seen before.

Firstly, we notice that while in the one dimensional case the concept of the derivative of a polynomial $Q\in \Q_{d,1}$ is straightforward and focuses our action inside the class $\Q_{d-1,1}$, in general dimension we need to speak about the gradient of $Q$, and thus the natural analogue of $\Q_{d-1,1}$ becomes
\beq\label{derivpolin}
\Q^{(1)}_{d,m}:=\left\{\nabla Q=\left(\frac{\partial Q}{\partial x_1},\,\ldots,\,\frac{\partial Q}{\partial x_m}\right)\,\big|\,Q\in \Q_{d,m}\right\}\:.
\eeq
Next, we introduce the analogue of \eqref{dqj}, that is for $I\subseteq\TT^m$ (dyadic) cube we set the  \emph{geometric factor} of the pair $(\nabla Q,I)$ as
\beq\label{dqjm}
\Delta_{\nabla Q}(I):=l(I)\,\||\nabla Q|\|_{L^{\infty}(I)}\:,
\eeq
where here $l(I)$ stands for the length of the cube $I$ while $|I|$ will remain the notation for the volume of $I$.

From here we have the obvious correspondence with the one dimensional algorithm for defining the family of tiles $\P$:

\textsf{Algorithm, general $m$.}

\begin{itemize}
\item fix a scale $k\in\N$ and consider the set $\I(k)$ of all dyadic cubes $I\subseteq \TT^m$ at scale $k$, \textit{i.e.} $l(I)=2^{-k}$.

\item fix $I\in\I(k)$ and consider a maximal separated set of $m-$tuple polynomials $\Q^{(1)}_{d,m}(I)\subset \Q^{(1)}_{d,m}$ such that
for any $\nabla Q_{1},\nabla Q_{2}\in \Q^{(1)}_{d,m}(I)$ one has
\beq\label{dqja}
\Delta_{\nabla Q_1-\nabla Q_2}(I)\geq 1\:.
\eeq

\item construct the uncertainty regions $\{P(\nabla Q,I)\}_{\nabla Q\in \Q^{(1)}_{d,m}(I)}$ \emph{partitioning} the set of polynomials $\Q^{(1)}_{d,m}$ such that
  \beq\label{dqjaa}
    P(\nabla Q,I)\subseteq\{\nabla\tilde{Q}\in \Q^{(1)}_{d,m}\,|\,\Delta_{\nabla\tilde{Q}-\nabla Q}(I)\leq \frac{1}{2}\}\:.
 \eeq

\item now we set a tile $P$ be a tuple of the form $(I, P(\nabla Q,I))$ where $I$ ranges through the set of all dyadic intervals inside the torus while $\nabla Q\in \Q^{(1)}_{d,m}(I)$. The set of all tiles is denoted with $\P$.
\end{itemize}

Now given $P=(I_P, P(\nabla Q_P,I_P))\in\P$, the analogue of \eqref{convbel} becomes
\beq\label{convbelm}
\nabla Q \in P\:\:\:\textrm{iff}\:\:\:\nabla Q\in P(\nabla Q_P,I_P))\,.
\eeq

Also, for $a>0$ we set
\beq\label{dqjaa}
    aP(\nabla Q,I)=\{\nabla\tilde{Q}\in \Q^{(1)}_{d,m}\,|\,\Delta_{\nabla\tilde{Q}-\nabla Q}(I)\leq \frac{a}{2}\}\:,
 \eeq
and denote with $aP$ the tuple  $(I_P, aP(\nabla Q_P,I_P))$. Extending naturally \eqref{convbelm}, we say
\beq\label{convbelma}
\nabla Q \in aP\:\:\:\textrm{iff}\:\:\:\nabla Q\in aP(\nabla Q_P,I_P))\,.
\eeq

Finally, the rest of the notations, remain as in our paper \cite{lvPC}, with the obvious modifications.

\section{\bf Discretization}\label{Discret}

This section presents the rather simple adaptations of the one dimensional dicretization process to our context.\footnote{For both historical lineage continuity (see \cite{c1} and \cite{f}) as well as argumentation clarity we present our proof on the $m-$ dimensional torus rather than on $\R^m$. However the latter situation follows similarly with no significant changes.} Recall the definition of the general Polynomial Carleson operator on the $m-$dimensional torus:
\beq\label{polcarlrecc}
C_{d,m}f(x):=\sup_{Q\in\Q_{d,m}}\left|
\,\int_{\TT^m}e^{i\,Q(y)}\,K(y)\,f(x-y)\,dy\,\right|\:.
\eeq
In what follows, for notational simplicity, we will refer to the operator $C_{d,m}$ as $T$.

We rewrite\footnote{For symmetry reasons we prefer to rewrite \eqref{polcarlrecc} in the equivalent form
$C_{d,m}f(x):=\sup_{Q\in\Q_{d,m}}\left|\,\int_{\TT^m}e^{i\,(Q(x)-Q(x-y))}\,K(y)\,f(x-y)\,dy\,\right|$.}
\beq \label{symgrm}
Tf(x)=\sup_{Q\in\Q_{d,m}}|T_{Q}f(x)|\:,
\eeq
with
\beq \label{TQ}
T_{Q}f(x)=\int_{\TT^m}{K(x-y)\:e^{i\,(Q(x)-Q(y))}\,f(y)\,dy}\,,
\eeq
and $Q\in\Q_{d,m}$ with $Q(x)=\sum_{|\b|\leq d}c_{\b}\,x^{\b}$ where
$\b=(\b_1,\ldots, \b_m)\in \N^m$ multi-index and as usual $x^{\b}=(x_1^{\b_1},\ldots, x_m^{\b_m})\in \TT^m$.

Now linearizing the supremum in $T$, we write
\beq \label{tq1}
Tf(x)=T_{Q_x}f(x)=\int_{\TT^m}K(x-y)\,e^{i\,(Q_x(x)-Q_x(y))}\:f(y)\,dy\:,
\eeq
where now $Q_x(y):=\sum_{|\b|\leq d}c_{\b}(x)\,y^{\b}$ with
$\{c_{\b}(\cdot)\}_{\b}$ measurable functions.

Further, we decompose our kernel $K$ as\footnote{For more on this decomposition one can consult \cite{sth}, Chapter 13.}
$$K(y)=\sum_{k\in\N} \psi_k(y)\:\:\:\:\:\:\:\:\:\forall\:\:y\in \TT^m,\,|y|\not=0\:,$$
where here $\psi_k(y):=2^{mk}\,\tilde{\psi}_k(2^k y)$ with each $\tilde{\psi}_k\in C^{1}_0(\R^m)$ being supported
in $\left\{u\in\R\:|\:2<|u|<8\right\}$ and satisfying uniformly in $k$
\beq \label{suppp}
 \|\partial_x^{\a}\tilde{\psi}_k \|_{L^{\infty}(\R^m)}\leq C\:\:\:\:\:\forall\:\:0\leq |\a|\leq 1,
\eeq
with $\int_{\R^m}\tilde{\psi}_k=0$ for all $k\in\N$.

Deduce that
\beq \label{tdec}
Tf(x)=\sum_{k\geq 0}T_{k}f(x):=\sum_{k\geq 0}\int_{\TT^m}e^{i\,(Q_x(x)-Q_x(y))}\,\psi_{k}(x-y)\,f(y)\,dy\:.
\eeq

Now for each $P=(I, P(\nabla Q,I))\in\P$ let
\beq \label{defEp}
E(P):=\left\{x\in I\:|\:\nabla_y Q_x (y)|_{y=x} \in P\right\}\,,
\eeq
and for $l(I)=2^{-k}$ ($k\geq0$), we
define the operators $ T_P$ on $L^2(\TT)$ by
\beq \label{deftp}
T_{P}f(x)=\left\{\int_{\TT^m}e^{i\,(Q_x(x)-Q_x(y))}\,\psi_{k}(x-y)\,f(y)\,dy\right\}\chi_{E(P)}(x)\:.
\eeq
With this, we deduce that
\beq \label{decT}
Tf(x)=\sum_{P\in\P}T_{P}f(x)\:.
\eeq

This ends our decomposition.
$\newline$

As in \cite{lvPC}, we end with this section with a useful remark that should be kept in mind for most of our later reasonings:

\begin{obs}\label{redu}
Taking $D$ to be the a very large positive integer depending polynomially on $m$ and $d$, and writing
$$\P=\bigcup_{j=0}^{D-1}\bigcup_{k\geq0}\P_{kD+j}\,,$$
we can assume from now on that the following \textit{scale separation condition} holds:
\beq \label{separat}
\eeq
if $P_j=(I_j, P_j(\nabla Q_{P_j},I_j))\in\P\:$ with
$j\in\left\{1,2\right\}$ such that
$|I_1|\not=|I_2|$ then either $|I_1|\leq 2^{-D}\:|I_2|$ or $|I_2|\leq 2^{-D}\:|I_1|$.
\end{obs}

\begin{comment}
\begin{obs}\label{redu}
We record here two facts that will be very useful in our later reasonings:
\begin{itemize}
\item For a tile $P=[\vec{\a},I_{P}]$, based on \eqref{suppp} and \eqref{deftp}, we deduce that
\beq \label{tptp*}
\operatorname{supp}\:T_P\subseteq I_P\:\:\:\:\:\textrm{and}\:\:\:\:\:\:\operatorname{supp}\:T_P^{*}\subseteq I_{P^*}\:,
\eeq
where here $T_P^{*}$ denotes the adjoint of $T_P$.

\item Taking $D$ to be the
smallest integer larger than $100d\log_{2}(100d)$ and splitting
$$\P=\bigcup_{j=0}^{D-1}\bigcup_{k\geq0}\P_{kD+j}\,,$$
we can assume from now on that the following \textit{scale separation condition} holds:
\beq \label{separat}
\eeq
if $P_j=[\vec{\a}_j,I_j]\in\P\:$ with
$j\in\left\{1,2\right\}$ such that
$|I_1|\not=|I_2|$ then either $|I_1|\leq 2^{-D}\:|I_2|$ or $|I_2|\leq 2^{-D}\:|I_1|$.
\end{itemize}
\end{obs}
\end{comment}

\section{\bf Quantifying the interactions between tiles}\label{Tileinteract}

Following the same lines as in \cite{lvPC} we adapt the estimates on the interaction between tiles to the general, higher dimensional case. The extensions are more or less straightforward with a natural correspondence in the one dimensional case. For transparency, our presentation mirrors the one in \cite{lvPC}.

\subsection{Properties of $T_{P}$ and $T_{P}^{*}$}
$\newline$

For $P=(I, P(\nabla Q_P,I))\in\P$ with $l(I)=2^{-k},\:k\in
\N$, we have \beq \label{v9}
\begin{array}{rl}
        &T_{P}f(x)=\left\{\int_{\TT^m}\:e^{i\,(Q_x(x)-Q_x(y))}\,\psi_{k}(x-y)\,f(y)\,dy\right\}\chi_{E(P)}(x)\:,  \\
    &T_{P}^{*}f(x)=\int_{\TT^m}\:e^{-i\,(Q_y(y)-Q_y(x))}\,\psi_{k}(y-x)\,\left(\chi_{E(P)}f\right)(y)\,dy\:.
\end{array}
\eeq
As expected, in direct correspondence with Observation 2 in \cite{lvPC}, we have the following heuristic:
\beq\label{loc}
\eeq
\begin{itemize}
\item the time-frequency
localization of $T_{P}$ is ``morally" given by $P$, i.e. should be essentially regarded as
$$\bigcap_{j=1}^m \{(x,\partial_j Q(x))|x\in I,\:\:Q\in\Q_{d,m},\:\: |\partial_j Q(x)-\partial_j Q_P(x)|\leq \frac{1}{2}\,l(I)^{-1}\}\:.$$
\item similarly, the time-frequency localization of $T_{P}^{*}$ is
``morally" given by $P^*$ with the analogue interpretation.
\end{itemize}

\subsection{Geometric factor of a tile relative to a polynomial} $\newline$

As before, given  $P=(I, P(\nabla Q_P,I))\in\P$ and $Q\in\Q_{d,m}$ one defines the geometric factor of $P$ relative to $Q$ (or $\nabla Q$) as
\beq\label{geomfact}
\left\lceil \Delta_{\nabla Q}(P)\right\rceil\:,
\eeq
where\footnote{Recall that given $x\in\R$ we let $ \left\lceil x\right\rceil:=\frac{1}{1+ |x|}$.}
\beq\label{v1}
\Delta_{\nabla Q}(P):=\inf_{\nabla Q_1\in P}\Delta_{\nabla Q-\nabla Q_1}(I)\:.
\eeq
Notice that we trivially have
$$\left\lceil\Delta_{\nabla Q}(P)\right\rceil\approx \frac{1}{1+l(I)\,\||\nabla Q-\nabla Q_P|\|_{L^{\infty}(I)}}\:.$$

\subsection{Van der Corput estimates for general dimension}

In this section we recall two results from \cite{sw} and very briefly discuss their adaptability to our context.

Throughout this section we assume we are given $Q\in\Q_{d,m}$ with $Q(x)=\sum_{|\b|\leq d}c_{\b}\,x^{\b}$. Let the \emph{size} of $Q$ be defined by
$$s(Q):=\sum_{|\b|\leq d}|c_{\b}|\,.$$
 With these, we have:
\begin{l1}(\cite{sw}) \label{Sww}
Let $\v\in C_0^{1}(\R^m)$ be a smooth function supported in the unit ball $B_m(0,1)$ and let $\Omega\subseteq B_m(0,1)$ be any convex set. Then, there exists a constant  $c=c(d,m)>0$ depending only on $d$ and $m$ such that
\beq\label{swest}
|\int_{\Omega} e^{i Q(x)}\,\v(x)\,dx|\leq c\, s(Q)^{-\frac{1}{d}}\,\|\v\|_{C^1(B_m(0,1))}\:.
\eeq
\end{l1}
One also has a good control over the size of the level sets:

\begin{l1}(\cite{sw}) \label{Sww1}
With the same notations as before, given any $\ep>0$, one has
\beq\label{swestt}
|\{x\in B_m(0,1)\,|\, |Q(x)|\leq \ep\}|\leq c\,\ep^{\frac{1}{d}}\,s(Q)^{-\frac{1}{d}}\:.
\eeq
\end{l1}

We now notice that both lemmas above are properly behaving under the action of dilation and translation symmetries. Moreover, we remark the following: given any $I\subseteq \TT^m$ (dyadic) cube and setting $P(y):=Q(l(I)y+c(I))$ one has the key relation
\beq\label{sizepol}
s(P)\gtrsim_{d,m}  \Delta_{\nabla P}\left([-\frac{1}{2},\frac{1}{2}]^m\right)=\Delta_{\nabla Q}(I)\:.
\eeq

From this and Lemma \ref{Sww} we immediately deduce

\begin{l1}\label{adaptpol}
Let $Q\in\Q_{d,m}$ and $I\subseteq \TT^m$ (dyadic) cube. Also, assume $\v_{I}\in C_0^{1}(10I)$ is a function adapted to $I$ with $\|\v_{I}\|_{L^{\infty}}\lesssim 1$ and $\Delta_{\nabla \v_I}(I):=l(I)\,\||\nabla \v_{I}|\|_{L^{\infty}(I)}\lesssim 1$. Then the following holds:
\beq\label{swestt1}
\left|\int_{\TT^m} e^{i Q(x)}\,\v_I(x)\,dx\right|\lesssim_{d,m}\, \left\lceil \Delta_{\nabla Q}(I)\right\rceil^{\frac{1}{d}}\,|I|\:.
\eeq
\end{l1}

\subsection{Control over the inner product} $\newline$

Consider $P_1=(I_1, P_1(\nabla Q_{P_1},I_1))\in\P$ and $P_2=(I_2, P_2(\nabla Q_{P_2},I_2))\in\P$. In what follows we quantify the output of the interaction
$$ |<T_{P_1}^{*}f, T_{P_2}^{*}g>|\,.$$
As expected, the output of the interaction will be controlled by ``the relative position" of
$P_1^{*}$ with respect to $P_2^{*}$ quantified in the definition below:

\begin{d0}\label{fact} [\textsf{Geometric factor associated to a pair of tiles}]

Let $P_1$ and $P_2$ be two tiles as above such that $I_{P_1}^{*}\cap I_{P_2}^*
\not=\es$. We define the {\bf geometric factor of the pair ($P_1,P_2$)} by
$$\left\lceil\Delta(P_1,P_2) \right\rceil\:,$$ where
$$\Delta(P_1,P_2)=\inf_{{\nabla Q_{1}\in P_1}\atop{\nabla Q_{2}\in P_2}} \Delta_{\nabla Q_1-\nabla Q_2}(\tilde{I}_{P_1}\cap \tilde{I}_{P_2})\:.$$
As before, we obviously have
$${\left\lceil\Delta(P_1,P_2)\right\rceil}\approx_{d,m}\max\left\{{\left\lceil\Delta_{\nabla Q_{P_1}}(P_2)\right\rceil},
\:{\left\lceil\Delta_{\nabla Q_{P_2}}(P_1)\right\rceil}\right\}\:.$$
\end{d0}

Now using Lemma \ref{adaptpol} above, one deduces the following:

\begin{l1}\label{interact}[\textsf{Tile interaction control}]

Let $P_1\:,\:P_2\:\in\P$. Then, with the above notations and conventions, we have
\beq\label{vdec}
|T_{P_1} T_{P_2}^{*}f(x)|\lesssim_{d,m}{\left\lceil
\Delta(P_1,P_2)\right\rceil}^\frac{1}{d}\:\frac{\int_{E(P_2)}|f|}{\max\left(|I_1|,|I_2|\right)}\,\chi_{E(P_1)}(x)\:. \eeq
\end{l1}

\begin{proof}
The proof is a straightforward application of Lemma \ref{adaptpol}, once one notices that
$$T_{P_1} T_{P_2}^{*}f(x) = \chi_{E(P_1)}(x)\,\int
(f\,\chi_{E(P_2)})(s)\:\K(x,s)\,ds\:,$$
where, assuming $l(I_1)=2^{-k_1}$ and $l(I_2)=2^{-k_2}$ with $k_1,\,k_2\in\N$, we let
$$\K(x,s)=\int e^{i\,
[(Q_s(s)-Q_s(y))- (Q_x(x)-Q_x(y))]}\;\psi_{k_1}(x-y)\:\varphi(y)\:\psi_{k_2}(y-s)\;dy\:.$$
\end{proof}

\begin{obs}\label{intlem}
In the one dimensional case, we got a more precise control over the $|<T_{P_1}^{*}f, T_{P_2}^{*}g>|$ expression in the form
\beq\label{v15} \left|\int
\tilde{\chi}_{I_{1,2}^c}T_{P_1}^{*}f\:\overline{T_{P_2}^{*}g}\:\right|\lesssim_{\:n,\:d,\:\ep_0}{\left\lceil
\Delta(P_1,P_2)\right\rceil}^{n}\:\frac{\int_{E(P_1)}
|f|\int_{E(P_2)}|g|}{\max\left(|I_1|,|I_2|\right)}\:\:\:\:\:\:\forall\:n\in
\N\:, \eeq
and
 \beq\label{v16} \int_{I_{1,2}}|
T_{P_1}^{*}f\:\overline{T_{P_2}^{*}g}|\lesssim_d{\left\lceil
\Delta(P_1,P_2)\right\rceil}^\frac{1-\ep_0}{d}\:\frac{\int_{E(P_1)}
|f|\int_{E(P_2)}|g|}{\max\left(|I_1|,|I_2|\right)}\:, \eeq
with $\tilde{\chi}_{I_{1,2}^c}$ a smooth variant of the corresponding
cut-off.

Here $I_{1,2}$ - called the ($\ep_0$)-critical intersection set - is essentially the region formed by the union of the maximal intervals $J\subseteq \tilde{I}_1\cap \tilde{I}_2$ for which one has $|\Delta_{\nabla Q_{P_1}-\nabla Q_{P_2}}(J)|\lesssim_d \Delta(P_1,P_2)^{\ep_0}$ where here $\ep_0\in (0,1)$. Based on Lemma \ref{Sww1}, one could obtain the higher dimensional analogue of \eqref{v15} and \eqref{v16}. However this would necessitate extra-technicalities that overweight the benefit - one does not need such an accurate description of the $ <T_{P_1}^{*}f, T_{P_2}^{*}g> $ interaction. Indeed, estimate \eqref{vdec} is more than enough for our final aim. For more on related considerations please see 1) in the Remarks section.
\end{obs}

\section{\bf The proof of the main theorem}\label{proofmainth}

No relevant modifications appear relative to the arguments presented in the corresponding section in \cite{lvPC}.

For convenience only, we will remind two of the key definitions involved later in our proof:

\begin{d0}\label{mass} [\textsf{Mass of a tile adapted to a given environment}]

Let $\A$ be a (finite) union of dyadic cubes in $\TT^m$ and $\p$ be a finite family of tiles. For  $P=(I, P(\nabla Q_P,I))\in\p$ with $I\subseteq\A$ we define the {\bf mass} of $P$ relative to the set of tiles $\p$ and the set $\A$ as being

\beq\label{v1} A_{\p,\A}(P):=\sup_{{P'= (I', P'(\nabla Q_{P'},I'))\in\:\p}\atop{I\subseteq
I'\subseteq \A}}\frac{|E(P')|}{|I'|}\:\left\lceil
\Delta(10P,\:10P')\right\rceil^{N} \eeq where $N$ is a fixed large
natural number.
\end{d0}

The qualitative concept that characterizes the
overlapping relation between tiles is given by

\begin{d0}\label{ord}[\textsf{Aiming for ``orderings"}]

   Let $P_j=(I_j, P_j(\nabla Q_{P_j},I_j))\in\P$ with $j\in\left\{1,2\right\}$. We say that
$\newline$ - $P_1\leq P_2$       iff       $\:\:\:I_1\subseteq I_2$
and      $\exists\:\:\nabla Q\in P_2$ such that $\nabla Q\in P_1\:,$ $\newline$ -
$P_1\trianglelefteq P_2$     iff     $\:\:\:I_1\subseteq I_2$  and
$\forall \:\:\nabla Q\in P_2$ we have $\nabla Q\in P_1\:.$

Also we say $P_1< P_2$ if $P_1\leq P_2$ and $|I_1|<|I_2|$. Similar statement for $\vartriangleleft$.
\end{d0}

\begin{obs}\label{ordmaxrel1}
Notice that $\leq$ is not an order relation while $\trianglelefteq$ it is. Also $P_1<P_2$ implies $2 P_1\vartriangleleft 2P_2$.
 \end{obs}

$\newline$
\subsection{\textbf{Partitioning $\P$}}\label{parttil}
$\newline$

Following with trivial adaptation to the higher dimensional case the stopping-time algorithm developed in \cite{lvPC},  we obtain a partition of our set of tiles into
\beq\label{decP}
\P=\bigcup_{n\in\N}\P_n\,,
\eeq
with each $\P_n$ being a set of tiles of mass $n$ relative to certain space regions. Our algorithm relies in a key fashion on
\begin{enumerate}
\item the concept of mass introduced in Definition \ref{mass};

\item a delicate analysis of  the level set of various counting functions involving the John-Nirenberg inequality.
\end{enumerate}

For concreteness we will only summarize the output of our algorithm. For this, following \cite{lvPC}, we simply quote (with the obvious adaptations) the following

\begin{d0}\label{dom}
 Let $\A=\bigcup \A_j$  and $\B=\bigcup_{l}\B_l$ be two collections of dyadic cubes inside $[0,1]^m$.
\begin{itemize}
 \item We say that
\beq\label{dom0}
\A\Subset\B\,,
\eeq
 iff given any two dyadic cubes $\A_j\in \A$ and $\B_l\in \B$ such that $\A_j\cap \B_l\not=\emptyset$ one has $\A_j\subseteq \B_l$ and also each $\A_j$ is contained in some $\B_l$.

\item We say that
\beq\label{dom1}
\A\prec\B\,,
\eeq
 iff each $\A_j$ is contained in some $\B_l$.

\item given an absolute constant $c>0$, we write
\beq\label{dom2}
\A\prec_{c}\B\,,
\eeq
 iff $\A\prec\B$ and for any $\B_l$ the following holds:
\beq\label{dom3}
|\bigcup_{\A_j\subseteq\B_l}\A_j|\leq e^{-c}\,|\B_l|\,.
\eeq
\end{itemize}
\end{d0}

$\newline$

\noindent\textbf{Output of the exceptional-set removing stopping-time algorithm}
$\newline$

\textit{There exists a collection of stopping-time dyadic cubes $\{\S_n\}_{n\in\N}$  and a corresponding collection of tiles $\{\P_n\}_{n\in\N}$ such that:
\begin{itemize}
\item each set $\S_n$ can be further partitioned as
\beq\label{sn}
\S_n:=\bigcup_{k\in\N} \S_{n}^k\,,
\eeq
with
\beq\label{sn1}
\S_{n}^{k+1}\Subset\S_{n}^k\,.
\eeq
\item moreover, for any $k'\leq k$
\beq\label{sn2}
\S_{n}^{k}\prec_{c (k-k')}\S_{n}^{k'}\,,
\eeq
and for any $n'\leq n$
\beq\label{sn3}
\S_{n}\prec\S_{n'}\,.
\eeq
\item one can partition the family of tiles
\beq\label{sn4}
\P=\bigcup_{n}\P_n\,,
\eeq
and further on partition each $\P_n$ as
\beq\label{sn4}
\P_n=\bigcup_{k}\P_n^k\,,
\eeq
such that
\begin{itemize}
\item for any $P\in \P_n^k$ we have $I_P\subseteq \S_{n}^k$ and $I_P\nsubseteq \S_{n}^{k+1}$;
\item for any $P\in \P_n^k$ we have $A_{\P,\S_n^k}(P)\in (2^{-n}, 2^{-n+1}]$;
\item there exists $c>0$ (large) absolute constant such that
\beq\label{sn5}
\|\sum_{P\in \p_{n}^{k,max}}\chi_{I_P}\|_{L^{\infty}}\leq c\,n\,2^n\,,
\eeq
where
\beq\label{maxn1}
\begin{array}{cc}
\p_{n}^{k,max}
:=\left\{P\in\P\,\big|\,\begin{array}{cc} P\:\textrm{maximal}\:\:\textrm{s.t.} \\I_P\subseteq \S_{n}^k\:\:\&\:\:\frac{|E(P)|}{|I|}> 2^{-n}
\end{array}\right\}\:.
\end{array}
\eeq
\end{itemize}
\end{itemize}}

The proof of this statement follows line by line the algorithm of partitioning the set of tiles $\P$ presented in Section 5.1.2 of \cite{lvPC}. We will thus not provide the details here but only mention that the sets $\{\S_{n}^k\}_{k\geq 1}$ are constructed inductively starting from $n=1$, $k=1$ as maximal disjoint collection of cubes formed from the analogue of the sets
$\{A_{n}^{l}[A_{n-1}^{j_{n-1}}\diamond A_{n-1}^{j_{n-1}+1}]\ldots[A_{1}^{j_1}\diamond  A_{1}^{j_1+1}]\}_{j_1,\ldots,j_{n-1},\,l}$ while $\P_n^k$ is formed in a similar fashion from the correspondent families $\{\p_n(A_{n}^{l}[A_{n-1}^{j_{n-1}}\diamond A_{n-1}^{j_{n-1}+1}]\ldots[A_{1}^{j_1}\diamond  A_{1}^{j_1+1}])\}_{j_1,\ldots,j_{n-1},\,l}$.

\subsection{\bf Main Proposition}
$\newline$
With this done, in direct correspondence with the similar statement in \cite{lvPC}, our main theorem follows by a simple application of triangle inequality paired with a geometric summation argument derived from
$\newline$

\noindent {\bf Main Proposition.}
{\it Fix $n\in\N$. Then there exist a constant $\eta=\eta(d,m)\in(0,\frac{1}{2})$ depending only
on $d$ and $m$ such that
$$\left\|T^{{\P}_n}f\right\|_{p}\lesssim_{p,d,m} \,2^{-n\,\eta(1-\frac{1}{p^*})}\,\left\|f\right\|_p\:,$$
for all $f\in L^p(\TT)$.}
$\newline$

\subsection{\bf Reduction of the main proposition}\label{Reductmainprop}

With the same notations and definitions from Section 6 in \cite{lvPC} at which we add Observation \ref{redu} in our paper one can follow the same reasonings  as in \cite{lvPC} (no significant modifications required) in order to reduce our Main Proposition above to the following statements:

$\newline$

\noindent\textbf{Proposition 1}\label{prop1} [\textsf{Control over a sparse forest}]

\textit{Let $\p\subseteq\P_n$ be a sparse forest. Then there exists $\eta=\eta(d,m)\in(0,1)$, depending only on the degree $d$ and dimension $m$,
such that for any $1<p<\infty$, we have
\beq\label{pro1}
\left\|T^{\p}\right\|_{p}\lesssim_{p,d,m} 2^{-n\,\eta\,(1-\frac{1}{p^*})}\:.
\eeq}
$\newline$

\noindent\textbf{Proposition 2}\label{prop2} [\textsf{Control over a (general) forest}]

\textit{Let $\p\subseteq\P_n$ be a forest. Then there exists $\eta=\eta(d,m)\in(0,1)$, depending only on the degree $d$ and dimension $m$,,
such that for any $1<p<\infty$ we have
\beq\label{pro2}
\left\|T^{\p}\right\|_{p}\lesssim_{p,d,m} 2^{-n\,\eta\,(1-\frac{1}{p^*})}\:.
\eeq}

\section{\bf The proofs of Propositions 1 and 2}\label{prop12}

The proof of Proposition 1 follows with no significant modifications the correspondent proof in \cite{lvPC}.

For the proof of Proposition 2 one can again follow the same reasonings as in \cite{lvPC} with the following adaptations corresponding to the lemmas dealing with
\begin{itemize}
\item the interaction of separated trees;

\item row-tree interaction.
\end{itemize}

In what follows we will only focus on these two modifications.

For this, we recall first several definitions:

\begin{d0}\label{sep}[\textsf{Separated trees}]

Fix a number $\d\in(0,1]$. Let $\p_1$ and $\p_2$ be two trees with
tops $P_j=(I_j, P_j(\nabla Q_{P_j},I_j))\in\P$ and $j\in\{1,2\}$. We say that $\p_1$ and
$\p_2$ are $\d^{-1}$-\textbf{separated} if either $I_1\cap
I_2=\es$ or else
\begin{itemize}
\item  $P=(I, P(\nabla Q_{P},I))\in\p_1\:\:\&\:\:I\subseteq
I_2\:\:\:\:\Rightarrow\:\:\:\left\lceil
\Delta(P,P_2)\right\rceil<\d\:,$
\item  $P=(I, P(\nabla Q_{P},I))\in\p_2\:\:\&\:\:I\subseteq
I_1\:\:\:\:\Rightarrow\:\:\:\left\lceil
\Delta(P,P_1)\right\rceil<\d\:.$
\end{itemize}
\end{d0}

\begin{l1}\label{sept}[\textsf{Interaction of separated trees}]

Let $\left\{\p_j\right\}_{j\in\left\{1,2\right\}}$  be two $\d^{-1}$-separated
trees with tops $P_j=(I_0, P_j(\nabla Q_{P_j},I_0))\in\P$. Then,
for any $f,\:g\in L^{2}(\TT)$, we have that
\beq\label{v21} \left|\left\langle
{T^{\p_1}}^*f,\,{T^{\p_2}}^*g\right\rangle\right|\lesssim_{m,d}{\d}^{\frac{1}{2d}}\left\|f\right\|_{L^{2}(\tilde{I}_0)}\left\|g\right\|_{L^{2}(\tilde{I}_0)}\:.
\eeq
\end{l1}
\begin{proof}

This is the analogue of Lemma 36 in \cite{lvPC}; one should notice the less refined form of our present statement. However this will be enough for our later estimates. For more on this, please see Observation \ref{septs1}.

We only sketch\footnote{For further details see the analogue proof of Lemma 2 in \cite{q} as well as that of Lemma 4 in \cite{f}.} the simple modifications of the argument presented in the proof of Lemma 36 in \cite{lvPC}:.

\begin{itemize}

\item Define a real-valued function $\f\in C_{0}^{\infty}(\R^m)$ with the following properties:
\begin{itemize}
\item $ supp\:\f\subset\left\{\frac{1}{4}\leq |x|\leq \frac{1}{2}\right\}$
\item $ \f\:is\: even $
\item $ |\hat{\f}(\xi)-1|\lesssim_{n} |\xi|^n \:\:\:\:\forall\:|\xi|\leq 1\:\:and\:n\: big\:enough $
\item $ |\hat{\f}(\xi)|\lesssim_{n} |\xi|^{-n}\:\:\:\:\forall\:|\xi|\geq 1$
    \end{itemize}

Next, define
$$d_{j}:=\min\{l(I_P)\,|\,P=(I_P, P(\nabla Q_{P},I_P))\in \p_j\}\:.$$

Now, for $j\in\left\{1,2\right\}$ let
\beq\label{fprop1}
\f_{j}(x)=(\d^{1/2}d_{j})^{-m}\f((\d^{1/2}d_{j})^{-1}x)\:.
\eeq
and define the operators
\beq\label{fprop2}
\tilde{\f}_{j}\::L^2(\R^m)\longrightarrow L^2(\R^m)\:\:\:\: by\:\:\tilde{\f}_{j}f=\f_{j}*f\,,
\eeq
and
\beq\label{fprop3}
\Phi_{j}\::L^2(\R^m)\longrightarrow L^2(\R^m)\:\:\:\:by\:\:\Phi_{j}= M_{\p_j} \tilde{\f}_{j} M^{*}_{\p_j}\:.
\eeq
where in the last line we define the general modulation operator
\beq\label{fprop03}
 M_{\p_j}:\:L^2(\R^m)\longrightarrow L^2(\R^m)\:\:\:\:by\:\:M_{\p_j}f(x)= e^{i\, Q_{P_j}(x)}\,f(x),,
\eeq
with $Q_{P_j}\in \Q_{d,m}$ being the ``antiderivative" central polynomial associated with the top tile $P_j=(I_0, P_j(\nabla Q_{P_j},I_0))$ and obeying $Q_{P_j}(0)=0$.

\item Following now similar reasoning with those in Lemma 2 in \cite{q}, for $j\in\left\{1,2\right\}$,
we decompose $T^{*}_{j}:={T^{\p_j}}^{*}$ as
\beq\label{dec}
T^{*}_{j}f=\Phi_{j}{T^{*}_{j}}f\:+\:\Omega_{j}f\,,
\eeq
and deduce that
\beq \label{om}
 \left\|\Omega_{j}\right\|_2\lesssim_{n}\d^{n}\,,
\eeq
while from an application of Lemma \ref{adaptpol} one gets
\beq \label{fi} \left|\left\langle\Phi_{1}{T^{*}_{1}}f
,\Phi_{2}{T^{*}_{2}}g\right\rangle\right|\lesssim_{m,d}\d^{\frac{1}{2d}}\left\|f\right\|_2\left\|g\right\|_2\:.
\eeq
finishing our proof.
\end{itemize}

\end{proof}

\begin{obs}\label{septs1}

In the one-dimensional we got a more refined version of the above lemma. Indeed, following Lemma 36 in \cite{lvPC}, one has that if $\left\{\p_j\right\}_{j\in\left\{1,2\right\}}$ are two $\d^{-1}$-separated
trees with tops $P_j=[\vec{\a}_j,I_0]$, then,
for any $f,\:g\in L^{2}(\TT)$ and $n\in \N$, one has that
\beq\label{v2111}
\left|\left\langle
{T^{\p_1}}^*f,\,{T^{\p_2}}^*g\right\rangle\right|\lesssim_{n,d}
{\d}^n\left\|f\right\|_{L^{2}(\tilde{I}_0)}\left\|g\right\|_{L^{2}(\tilde{I}_0)}+
\left\|\chi_{I[c]}{T^{\p_1}}^*f\right\|_2\left\|\chi_{I[c]}{T^{\p_2}}^*g\right\|_2\:,
\eeq
where here the component $I[c]$ - called the  \emph{critical
intersection set} is the one responsible for the correspondent weak decay form appearing in \eqref{v21}. As we will see later, we don't actually need anything stronger than \eqref{v21}.
\end{obs}

\begin{d0}\label{nor}[\textsf{Normal tree}]

A tree $\p$ with top $P_0=(I_0, P_j(\nabla Q_{P_0},I_0))$ is called \textbf{normal} if for any
$P=(I, P(\nabla Q_{P},I))\in\p$ we have $100 I\cap (I_0)^c=\emptyset\;.$
\end{d0}

\begin{obs}\label{norob} Notice that if $\p$ is a normal tree as above
then $$supp\:{{T^{\p}}^{*}f}\subseteq I_0\:.$$
\end{obs}

\begin{d0}\label{row} [\textsf{Row}]

A \emph{row} is a collection $\p=\bigcup_{j\in \N}\p^{j}$ of normal
trees $\p^{j}$ with tops
$P^{j}_0=(I_0^j, P_j(\nabla Q_{P_0^j},I_0^j))$ such that the dyadic cubes
$\left\{I^j_0\right\}$ are pairwise disjoint.
\end{d0}

\begin{l1}\label{rt}[\textsf{Row-tree interaction}]

Let $\p$ be a row as above, let $\p'$ be a tree with top $P'_0=(I'_0, P_j(\nabla Q_{P'_0},I'_0))$ and suppose
that $\forall\: j\in \N$, $I_0^{j}\subseteq I'_0$ and $\p^{j},\:\p'$
are $\d^{-1}$separated trees.

Then for any $f,\:g\in L^{2}(\TT)$ we have that
\beq \label{rotre}
\left|\left\langle {T^{\p'}}^*f,{T^{\p}}^*g\right\rangle\right|
\lesssim_{m,d}{\d}^{\frac{1}{2d}}\,\left\|f\right\|_{2}\left\|g\right\|_{2}\:.
\eeq
\end{l1}

\begin{proof}
The proof is straightforward and follows from Lemma \ref{sept}, Observation \ref{norob} and Definitions \ref{nor} and \ref{row}. The interested reader might want to check for comparison the proof of the corresponding Lemma 40 in \cite{lvPC}.
\end{proof}

We now state the Main Lemma whose statement and proof follows line by line the corresponding ones in \cite{lvPC}.
$\newline$

\noindent\textbf{Main Lemma.} {\it Let $\p\subset \P_n$ be an $L^{\infty}-$forest of generation $n$.

\noindent Then there exists $\eta=\eta(d,m)\in (0,1)$ such that
\beq \label{mle}
\|T^{\p}f\|_2\lesssim_{d,m} 2^{-\frac{n}{2}\,\eta}\,\|f\|_2\:.
\eeq}

\begin{proof}
Below, we only present the outline of the proof of this result - for more details please check the proof of the Main Lemma in \cite{lvPC}:

\begin{itemize}
\item based on the definition of an $L^{\infty}-$forest we can decompose $\p$ as
\beq \label{decp}
 \p=\bigcup_{j=1}^{c\,2^n} \r_j\,,
\eeq
 with each $\r_j$ being a collection of pairwise spatially disjoint trees.

\item  on each of the trees belonging to a given $\r_j$ we perform two operations: 1) we trim the bottom - the first $100\,n\,d$  minimal tiles; 2) we remove its boundary component. The total excised collection of tiles (over all the trees in $\p$) may be decomposed in at most $c\,n$ sets of incomparable tiles to which one applies Proposition 1.

\item  with this we reduce our discussion to the situation in which \eqref{decp} represents the row-decomposition of an $2^{100\,n\,d}$-separated $L^{\infty}-$forest $\p$; reached at this point we notice that the operators $\{T^{\r_j}\}_j$ are almost orthogonal.

\item thus, for $k\not=j$, we have
\begin{itemize}
\item $\|{T^{\r_k}}^{*}\,T^{\r_j}\|_{2\mapsto 2}=0$ - as a consequence of the pairwise disjointness of the sets $\{\textrm{supp}\,T^{\r_j}\}_j$;

\item $\|T^{\r_k}\,{T^{\r_j}}^{*}\|_{2\mapsto 2}\lesssim 2^{-5\,n}$ - based on the ($2^{100\,n\,d}$)-separateness hypothesis to which we apply Lemma \ref{rt}. Here is the key point where we notice that one does not need anything more sophisticated than the form given by \eqref{rotre}.
\end{itemize}

\item Using now the last item and applying the Cotlar-Stein lemma and the single tree estimate we conclude that \eqref{mle} holds.

\end{itemize}
\end{proof}

With these done, one can apply the same reasonings as in \cite{lvPC} in order to deduce that Proposition 2 holds.

\section{\bf Remarks}\label{remark}

1) The first remark concerns the more intricate nature of the exposition of the one dimensional versus the higher dimensional case. Many of the arguments in \cite{lvPC} could be written in a more succinct form. However, the author's intentions in \cite{lvPC} were to provide an approach that is detailed, self-contained and, most importantly, offers a transparent antithesis between the situations in which standard wave-packet analysis techniques apply versus situations in which one encounters a new manifestation that is specific to the generalized wave-packet setting. With respect to the latter item, we exemplify with the statement of, say, Lemma 8 in \cite{lvPC}. There, one can see a dichotomy between the fast decay estimate obtained in the regions far away from the ``intersection of the geometric tiles" (similar with the original Fefferman setting, see \cite{f}) and the minimal decay obtained over the regions where one has overlapping of the tiles, situation that is specific only to the higher than one degree polynomial phases setting. In order to unravel this antithesis one has to perform a fine analysis\footnote{Including level set estimates, min/max properties, growth of derivatives etc.} of the behavior of a polynomial on a given interval - see Lemma 3 in \cite{lvPC}.

In contrast with this, the analogue result of  Lemma 8 in \cite{lvPC}, that is Lemma \ref{interact} here, focuses strictly on some (minimal) decay of the interaction. As a consequence of the rougher estimates one can drop the more intricate analysis analogue to Lemma 3 in \cite{lvPC} and simply quote a Van der Corput type estimate in the form of Lemma \ref{adaptpol} here. Similar situations appear when one compares the statements and proofs of Lemmas 36 and 40 in \cite{lvPC} versus their present analogue Lemmas \ref{sept} and \ref{rt} respectively.

However, passing over the format of the presentation, the essence of the approach remains the same and follows the fundamental ideas introduced in \cite{q} and \cite{lvPC}:

\begin{enumerate}
\item the definition of the uncertainty regions/tiles based on the \emph{relational approach} defined in \cite{q} (see also the Introduction in \cite{lvPC} and the reference to \cite{Fefunc});

\item the stopping time algorithm on removing the exceptional sets followed by the associated tile partitioning $\P=\bigcup_n\P_n$ adapted to a specialized concept of mass of a tile;

\item the decomposition of each of the families $\P_n$ into a controlled number of forests.
\end{enumerate}

Indeed, these three elements are all present in the higher dimensional case - see both \cite{zk} and the current paper - following with simple modifications/adaptations the reasonings in \cite{lvPC}.

2) Our intention here is to make transparent the connections between some of the few seemingly different approaches that are developed in \cite{lvPC} and our present work as opposed to the work in \cite{zk}.

2.1) Our first comment refers to the definition of tiles. In our work the key quantity involved in the tiling of the time-frequency domain is given by the geometric factor of the pair\footnote{Recall that in our context $Q\in \Q_{d,m}$ is a polynomial in $\R^m$ of degree at most $d$ while $I\subset \TT^m$ is a dyadic cube.} $(\nabla Q,I)$ defined in \eqref{dqjm} and recalled below:
\beq\label{recdqjm}
\Delta_{\nabla Q}(I):=l(I)\,\||\nabla Q|\|_{L^{\infty}(I)}\:.
\eeq
In \cite{zk}, the author chooses to work with a pair of the form\footnote{Same notations as above.} $(Q, I)$ and the quantity denoted by $\|Q\|_{I}$ and defined as
\beq\label{zk1}
\|Q\|_{I}:=\sup_{x,\,x'\in I} |Q(x)-Q(x')|\:.
\eeq
With this, the tile partitioning in \cite{zk} relative to the spacial location $I$ revolves around the idea of maximal $1-$separated (central) polynomials $Q\in \Q_{d,m}$ relative to the quantity defined in \eqref{zk1}. This type of requirement is in fact morally equivalent with our condition \eqref{dqja}. Indeed as a support for our claim we record the following simple representation formula:
\beq\label{obgr}
\eeq
$$Q(x)-Q(x')=$$
$$\int_{x_1'}^{x_1}\frac{\partial Q}{\partial y_1}(y_1,x_2',\ldots, x_m') dy_1\,+\,\ldots\,+\int_{x_m'}^{x_m}\frac{\partial Q}{\partial y_m}(x_1,x_2,\ldots, x_{m-1}, y_m) dy_m\,,$$
where as expected, for a generic $x\in\R^m$ we have denoted $x=(x_1,\ldots, x_m)$.

With this one notices the implication
\beq\label{imp}
\Delta_{\nabla Q-\nabla Q_P}(I)\lesssim 1\:\:\:\Rightarrow\:\:\:\|Q-Q_P\|_{I}\lesssim 1\:,
 \eeq
that realizes the translation between the tile discretization here and the one in \cite{zk}.

2.2) The definition of an antichain in \cite{zk} corresponds in fact to the definition of a family of incomparable tiles - see Definition 14 in \cite{lvPC}.

2.3) Lemma 3.8 in \cite{zk} is a reorganized and very compressed form of the stopping time algorithm on removing the exceptional sets introduced in Section 5.1.2. in \cite{lvPC}. We embrace this succinct form presentation in our present paper too - though with several key modifications - see the ``Output of the exceptional-set removing stopping-time algorithm" in Section 5.1.

2.4) Definition 3.18  in \cite{zk} of what the author there calls ``a Fefferman forest of level n and generation k" is essentially - up to a suitable log loss - an example of what we called an $L^{\infty}-$forest - see Definition 21 in \cite{lvPC}.

2.5) Proposition 3.23 in \cite{zk} follows with minor changes/language adaptations the argument in Section 6.2 in \cite{lvPC}.

2.6) Lemma 4.5 and Proposition 4.8 in \cite{zk} are treated both as part of Proposition 1 in \cite{lvPC}.

2.7) Lemmas 5.6 and 5.12 in \cite{zk} correspond to Lemma 36 and Lemma 45 respectively in \cite{lvPC}.

2.8) Lemma 5.13 in \cite{zk} corresponds to the proof of Proposition 2 - the $L^2$ bound case - see Section 7.3.1. in \cite{lvPC}.

2.9) We end this remark by noticing the interesting form of Lemma 4.1 in \cite{zk} which unravels in a beautiful manner the connection between the cancelation encoded in the oscillatory integral term $\int e^{i Q}\,\psi$  and the degree of smoothness encoded in the supremum over the variation of the function $\psi$ over the set of intervals whose size relates with $\|Q\|_I$ where $I=\textrm{supp}\,\psi$. This explains how the regularity at the critical scale of the function $\psi$ transfers into the decay of the initial expression  $|\int e^{i Q}\,\psi|$.

3) In this remark we want to explain a subtle point in defining the mass of a tile that affects the correctness of part of the reasonings in \cite{zk}. To be more precise, in contrast with our Definition \ref{mass} here, in \cite{zk} the author defines the mass of a tile
$P\in\P$ as simply given by
\beq\label{mt}
A_{s}(P):=\sup_{\l\geq 1}\frac{E(\l P)}{|I_P|}\,\l^{-N}\:,
\eeq
for some suitable fix natural number $N\in\N$.

While tempting due to the simplifications that it would have brought to the tile partitioning $\P=\bigcup_{n\in\N}\P_n$ and implicitly to the exceptional-set removing stopping-time algorithm, definition \eqref{mt} has a major flaw: if defined in this way, the mass of a tile $P$ is \emph{insensitive to how the tile $P$ embeds in the ambient universe}, specifically, there is no direct relationship between the mass of $P$ and the mass of tiles living at \emph{different} scales. In particular, one looses the original \emph{monotonicity} property of the mass, \textit{i.e.}
\beq\label{mt1}
P\leq P'\:\:\nRightarrow\:\:A_{s}(P)\geq A_{s}(P') \:.
\eeq

This aspect has in turn a major impact on the consistency of the tree selection algorithm:

Given $\p\subset\P$ a tree and adopting for a mass of a tile definition \eqref{mt}, one looses any control over the convexity of a tree with uniform mass; more precisely,\footnote{For simplicity one may want to imagine that the entire discussion takes place for $m=d=1$, hence simply for the standard Fefferman case.} defining
$$\p_n(x):=\{P\in\p\,|\,x\in I_P\:\:\:\&\:\:\:A_{s}(P)\in [2^{-n-1},\,2^{-n})\}\,,$$
where here $x\in I_0$ with $P_0$ the top of $\p$, we notice that the set of scales associated with the spacial intervals of the tiles in $\p_n(x)$ can have arbitrarily many gaps each one of arbitrary size. Indeed, if one departs with $\p$ being - say the time-frequency portrait of the standard Hilbert transform - then denoting with $$A_n(x,\p):=\{k\,|\,\exists\:P\in\p_n\:\:s.t\:\:l(I_P)=2^{-k}\:\:\&\:\:x\in\cdot I_P\}\,,$$ one has that $A_n(x,\p)$ can essentially be any subset of $\N$.

One could in principle try to overcome this difficulty by developing a theory for non-convex trees. The observation about an alternative approach that uses the concept of \emph{generalized trees} - that is trees that are not required to also be convex - was explicitly made in \cite{lvPC}, see remark 3 in Section 8. However, in that context our intention was to explain that allowing to work with generalized trees instead of standard trees one could obtain a simplification of the  exceptional-set removing stopping-time algorithm introduced in Section 5.1. of  \cite{lvPC}, but, \emph{importantly},  still using the ``smoothing" definition of the mass as provided by Definition \ref{mass}.

If instead, one intends to work with the mass definition \eqref{mt}, one needs to take in account the fact that the structure of the generalized trees of uniform mass becomes much rougher. Thus, a correct adaptation of the exceptional-set removing stopping-time algorithm and of the construction of the families $\P_n$ in \cite{lvPC} to the usage of \eqref{mt} would necessitate a highly nontrivial amount of technicalities and work that is absent in \cite{zk}. Moreover, even if successful, such a strategy would be unnecessarily complicated and essentially reworking - even if in a disguised manner - the exceptional-set removing stopping-time algorithm based on the concept of mass given by Definition \ref{mass}.


\begin{thebibliography}{10}


\bibitem{c1}
Lennart Carleson.
\newblock On convergence and growth of partial sumas of {F}ourier series.
\newblock {\em Acta Math.}, 116:135--157, 1966.


\bibitem{f}
Charles Fefferman.
\newblock Pointwise convergence of {F}ourier series.
\newblock {\em Ann. of Math. (2)}, 98:551--571, 1973.

\bibitem{Fefunc}
Charles~L. Fefferman.
\newblock The uncertainty principle.
\newblock {\em Bull. Amer. Math. Soc. (N.S.)}, 9(2):129--206, 1983.


\bibitem{hu}
Richard~A. Hunt.
\newblock On the convergence of {F}ourier series.
\newblock In {\em Orthogonal {E}xpansions and their {C}ontinuous {A}nalogues
  ({P}roc. {C}onf., {E}dwardsville, {I}ll., 1967)}, pages 235--255. Southern
  Illinois Univ. Press, Carbondale, Ill., 1968.


\bibitem{lt3}
Michael Lacey and Christoph Thiele.
\newblock A proof of boundedness of the {C}arleson operator.
\newblock {\em Math. Res. Lett.}, 7(4):361--370, 2000.


\bibitem{q}
Victor Lie.
\newblock The (weak-{$L^2$}) boundedness of the quadratic {C}arleson operator.
\newblock {\em Geom. Funct. Anal.}, 19(2):457--497, 2009.

\bibitem{lvPC}
Victor Lie.
\newblock The {P}olynomial {C}arleson operator.
\newblock {\em Arxiv: https://arxiv.org/abs/1105.4504., submitted}.


\bibitem{s2}
Elias~M. Stein.
\newblock Oscillatory integrals related to {R}adon-like transforms.
\newblock In {\em Proceedings of the {C}onference in {H}onor of {J}ean-{P}ierre
  {K}ahane ({O}rsay, 1993)}, number Special Issue, pages 535--551, 1995.


\bibitem{sth}
Elias~M. Stein.
\newblock Harmonic Analysis: real variable methods, orthogonality, and oscillatory integrals.
\newblock With the assistance of Timothy S. Murphy. {\em Princeton Mathematical Series, 43. {M}onographs in {H}armonic {A}nalysis, III. {P}rinceton {U}niversity {P}ress, {P}rinceton, NJ, 1993}.


\bibitem{sw}
Elias~M. Stein and Stephen Wainger.
\newblock Oscillatory integrals related to {C}arleson's theorem.
\newblock {\em Math. Res. Lett.}, 8(5-6):789--800, 2001.


\bibitem{zk}
Pavel Zorin-Kranich.
\newblock {M}aximal polynomial modulations of singular integrals.
\newblock {\em Arxiv: https://arxiv.org/abs/1711.03524}.


\end{thebibliography}
\end{document}